\documentclass[11pt]{article}

\usepackage[verbose=true,letterpaper]{geometry}

\usepackage[utf8]{inputenc} 
\usepackage[T1]{fontenc}    
\usepackage{hyperref}       
\usepackage{url}            
\usepackage{booktabs}       
\usepackage{amsfonts}       
\usepackage{nicefrac}       
\usepackage{microtype}      
\usepackage{lipsum}
\usepackage{fancyhdr}       
\usepackage{graphicx}       
\usepackage{amsmath}
\usepackage{amsthm}
\usepackage{xcolor}
\usepackage{comment}
\usepackage{setspace}
\usepackage{cleveref}

\graphicspath{{media/}}     

\theoremstyle{definition}

\author{Matan Eliashar and Nati Linial}
\date{}
\title{A Note on Fermat's Last Theorem for $n=4$}

\begin{document}
\maketitle

Fermat’s last theorem (FLT) famously states that the equation $x^n+y^n=z^n$ has no solution
in positive integers $x, y, z$ for any integer exponent $n>2$. But does this theorem have
a quantitative version? Specifically, for a given $n\ge 3$, how large can $d_n(\cdot)$ be if it is true
that $|z^n-x^n-y^n| \ge d_n(z)$ for every three positive integers $z> x,y>0$?

Exponent $n=4$ is the easiest case of FLT and was already known to Fermat. In fact, his proof shows that
even the equation $x^4+y^4=z^2$ has no solution in positive integers $x, y, z$, where $z> x^2, y^2$ \cite{hardy75}.
We consider here approximate solutions of this equation.
A computer search has yielded the following triplets:
\begin{center}
\begin{tabular}{ | c | c | c | c | c | } 
\hline
 n & 0 & 1 & 2 & 3 \\ 
\hline
 x & 22 & 1058 & 50806 & 2439746 \\ 
\hline
 y & 23 & 1103 & 52967 & 2543519 \\
\hline
 z & 717 & 1653213 & 3812308653 & 8791182100413 \\
\hline
\end{tabular}
\end{center}

with
\begin{equation}\label{eqn:general}
x_n^4 + y_n^4 - 8 = z_n^2
\end{equation}
These turn out to be the initial terms of an infinite sequence
of triplets $\left(x_n, y_n, z_n\right)_{n=1}^\infty$, that satisfy Equation (\ref{eqn:general}).
Here $x_n, y_n$ are defined by the above initial conditions and the recurrences
\begin{equation}\label{eqn:xy}
x_n = 48 \cdot x_{n-1} + x_{n-2}~~;~~y_n = 48 \cdot y_{n-1} + y_{n-2}
\end{equation}
whereas
\begin{equation}\label{eqn:z}
z_n=2306 \cdot z_{n-1} - z_{n-2} + (-1)^n \cdot 192
\end{equation}
These recurrences and the initial conditions clearly imply that $x_n, y_n, z_n$ 
are always positive integers. Using the three-term recurrence (\ref{eqn:xy}) and the initial conditions,
we get the closed-form expressions for $x_n, y_n$: 

\begin{equation}\label{eqn:explicit:xy}
x_n=a\lambda_1^n+b\lambda_2^n~~;~~y_n=c\lambda_1^n+d\lambda_2^n
\end{equation}
where \[\lambda_1=\left( 24 + \sqrt{577} \right) ~,~ \lambda_2=\left( 24 - \sqrt{577} \right)\]
are the two roots of the quadratic $\lambda^2=48\lambda+1$ corresponding to the recurrence (\ref{eqn:xy}), and
\[a=\left( 11 + \frac{265}{\sqrt{577}} \right) \; , \; b=\left( 11 - \frac{265}{\sqrt{577}} \right)  \; , \; c= \left( \frac{23}{2} + \frac{551}{2\sqrt{577}} \right) \; , \; d = \left( \frac{23}{2} - \frac{551}{2\sqrt{577}} \right)\]
Note also that $\lambda_1 \cdot \lambda_2 = -1$, the constant term in that quadratic.\\ 
Consequently:
\[x_n^4=a^4\lambda_1^{4n} + \left( -1 \right)^n 4a^3 b \lambda_1^{2n} + 6a^2b^2 + \left( -1 \right)^n 4ab^3\lambda_1^{-2n} + b^4\lambda_1^{-4n} \]
and hence:
\[ x_n^4 + y_n^4 - 8 =\]
\[ \left( a^4 + c^4 \right) \lambda_1^{4n} + \left( -1 \right)^n \left( 4a^3 b + 4c^3 d\right) \lambda_1^{2n} + \left( 6a^2b^2 + 6c^2d^2 - 8\right) + \left( -1 \right)^n \left( 4ab^3  + 4cd^3 \right) \lambda_1^{-2n} + \left( b^4 + d^4 \right) \lambda_1^{-4n} \]
The recurrence (\ref{eqn:z})
that $z_n$ satisfies is a bit less standard, but its closed-form solution is not hard to derive. Namely,
\begin{equation}\label{eqn:explicit:z}
z_n=e\mu_1^n+f\mu_2^n + (-1)^n \cdot g
\end{equation}
where \[\mu_1=\left( 1153 + 48\cdot\sqrt{577} \right) ~,~ \mu_2=\left( 1153 - 48\cdot\sqrt{577} \right)\]
are the two roots of the quadratic $\mu^2=2306\mu-1$ corresponding to three-term part of recurrence (\ref{eqn:z}), 
and
\[e=\left( \frac{413661}{1154}+\frac{17221}{1154}\sqrt{577} \right) \; , \; f=\left( \frac{413661}{1154}-\frac{17221}{1154}\sqrt{577} \right) \; , \; g=\left( \frac{48}{577} \right).\]
Note that  $\mu_1 = \lambda_1^2 ~,~ \mu_2 = \lambda_2^2 = \lambda_1^{-2}$,
whence also $\mu_1 \cdot \mu_2 = 1$.
Consequently:
\[z_n^2= e^2 \lambda_1^{4n} + 2eg\left( -1 \right)^n \lambda_1^{2n} + \left( 2ef+g^2\right) +  2fg\left( -1 \right)^n \lambda_1^{-2n} + f^2 \lambda_1^{-4n}\] 
\begin{flushleft}
Finally, we verify the following five equalities:
\end{flushleft}
\begin{itemize}
\item 
$e^2=a^4+c^4$
\item
$f^2=b^4+d^4$
\item
$2eg=4a^3b+4c^3d$
\item
$2fg=4ab^3+4cd^3$
\item
$2ef + g^2=6a^2b^2+6c^2d^2-8$
\end{itemize}
which together imply Equation (\ref{eqn:general}).
\newline
\bibliographystyle{plain}
\bibliography{bibliography.bib}
\end{document}